\documentclass {article}

\usepackage{amsfonts,amsmath,latexsym}
\usepackage{amstext}
\usepackage {epsf}
\usepackage {epsfig}
\usepackage[scanall]{psfrag}
\usepackage {graphicx}

\newcommand {\R}{\mathbb{R}}
\newcommand {\Ss}{\mathbb{S}}
\newcommand {\C}{\mathbb{C}}

\newcommand {\del}{\partial}

\newcommand {\sech}{\operatorname{sech}}
\newcommand {\dist}{\operatorname {dist}}
\newcommand {\Span}{\operatorname {span}}
\newcommand {\imaginary}{\operatorname {\Im m}}
\newcommand {\real}{\operatorname {\Re e}}
\newcommand {\calC}{\mathcal{C}}
\newcommand {\calD}{\mathcal {D}}
\newcommand {\calT}{\mathcal{T}}
\newcommand {\calG}{\mathcal{G}}

\newcommand {\calL}{\mathcal{L}}
\newcommand {\calV}{\mathcal{V}}

\newcommand {\Lap}{\operatorname{\Delta}\nolimits}
\newcommand{\Lip}{\operatorname{Lip}\nolimits}

\newtheorem {thm} {Theorem}

\newtheorem {lemma} [thm] {Lemma}

\newtheorem {rmk} {Remark}

\begin {document}
\title {A Capture Problem in Brownian Motion and Eigenvalues of Spherical
Domains}
\author {Jesse Ratzkin \\ University of Connecticut \\ ratzkin@math.uconn.edu
  \and Andrejs Treibergs \\ University of Utah \\ treiberg@math.utah.edu}
\maketitle
\begin {abstract}
We resolve a question of Bramson and Griffeath by showing
that the expected capture time
of four independent Brownian predators pursuing one Brownian prey
on a line is finite. Our main tool is an
eigenvalue estimate for a particular spherical domain, which we obtain
by a coning construction and domain perturbation.
\end {abstract}

In this paper, we examine the expected capture time of a single Brownian
prey pursued by $n$ independent Brownian
predators. All motion is restricted to a line. Bramson and Griffeath
\cite {BG} first considered this
problem, and estimated the capture time in various circumstances. In
particular, they showed that if at
time $t=0$ there are predators on both sides of the prey then the expected
capture time is finite. For
this reason, we will assume that the initial position of the prey is
$x_0(0) =
1$ and the initial positions
of the predators are $x_1(0)= \dots = x_n(0) = 0$. In this case,
\cite {BG} showed that the
expected capture time is infinite for $n = 1,2,3$, and conjectured that it is
finite for $n \geq 4$ (as
indicated by simulations). Li and Shao \cite{LS} showed
that the expected capture time is finite for $n \geq 5$.
Using a similar method, we resolve the remaining case by showing that the
expected capture time is finite for $n=4$.
\begin {thm} \label {main-thm}
Let $x_0(t)$ be pursued by $x_1(t), \dots, x_n(t)$, where $x_0, \dots,
x_n$ are independent, standard Brownian motions on a line, $x_0(0) = 1$, and
$x_1(0) = \cdots x_n(0) = 0$. Then the expected capture time is finite
if and only if $n \geq 4$.
\end {thm}
One key difference between Li and Shao's method in \cite{LS} and our
method here is that they consider a difference process, while we
do a geometric splitting in spherical polar coordinates  to reduce 
the dimension of the problem. See Section
\ref{separate-sec} for more details on our dimension reduction.

One can reformulate the capture problem described above
in terms of exit of a Brownian particle in $\R^{n+1}$ from a specific cone. We
denote the position of the prey at time
$t$ as $x_0(t)$ and the position of the $j$th predator at time $t$ as
$x_j(t)$. By our choice of initial
conditions, the initial position of the Brownian particle $x(t) = (x_0(t),
x_1(t), \dots, x_n(t))$ is
$x(0) = (1,0,\dots,0)$. The event of capture is then equivalent to the
Brownian particle $x(t)$ first leaving
the cone
$$
\calC_{n+1} := \{(x_0, x_1,\dots, x_n) \mbox{ } | \mbox{ }
x_0 \geq x_j, j = 1, \dots, n \},
$$
with $x(0) = (1,0,\dots, 0)$, and so we must estimate the expected
first exit time 
of a Brownian particle from the cone $\calC_{n+1}$, with the starting
position $(1,0,\dots,0)$.

DeBlassie \cite {DB} developed the theory of estimating exit times for
Brownian motion from cones in Euclidean space. Let $C = \{(r,\theta)
\mbox{ }|\mbox{ } r \geq 0, \theta \in D \subset \Ss^n\}$ be the
cone over a domain $D \subset \Ss^n$. Also let $\tau_x$ be the
exit time from $C$ of a Brownian particle with starting position $x$,
and let
$\mathbb{P}(\tau_x >t)$ be the probability that $\tau_x > t$. DeBlassie
showed $\mathbb{P}(\tau_x > t) \sim c(x) t^{-a(n)}$, where
\begin {equation} \label{decay-exp}
2a(n) = \left [ \left (\frac{n-2}{2} \right )^2 + \lambda_1(D)
\right ]^{1/2} - \frac{n-1}{2}.
\end{equation}
Here $\lambda_1(D)$ is the first Dirichlet eigenvalue of $D$. In the
particular case we are interested in, the expected exit time of a Brownian
particle from $\calC_{n+1}$ is finite if and only if $a(n) > 1$,
which reduces to
\begin {equation} \label{eigen-est1}
\lambda_1(\calD_n) > 2n + 2,
\end {equation}
where $\calD_n = \calC_{n+1} \cap \Ss^n$.

Our method for proving Theorem \ref{main-thm} is to estimate the first
eigenvalue
of $\calD_n$ using the monotonicity property of eigenvalues, a coning
construction, and domain perturbation. The rest of the paper proceeds as
follows. In Section
\ref{geometry-sec} we discuss the geometry of $\calC_{n+1}$, $\calD_n$, and
related regions. Section \ref{separate-sec} contains the separation of
variables
background one needs
to estimate the expected capture time. For the reader's convenience, we also
include the proof that the expected capture time is infinite for $n = 1,2,3$
predators in this section. We prove the $n=4$ eigenvalue estimate
in Section \ref{eigen-est-sec}. Finally, we describe a numerical computation of
$\lambda_1(\calD_3)$ and a lower bound for $\lambda_1(\calD_4)$
in Section~\ref{numer-sec}.

{\scriptsize {\sc Acknowledgement:} Davar Khoshnevisan told us about this
problem. We thank him and Pedro Mendez for lending  their ears
and expertise in
this project. A. T. thanks Frank Stenger for suggesting the algorithm of
Section~\ref{numer-sec}.}

\section {Geometry of the cone $\calC$} \label{geometry-sec}

The cone $\calC_{n+1}$ and its spherical angle $\calD_n$ have much symmetry.
First observe that $\calC_{n+1}$ contains the line $\calL$ spanned by
$(1,1,\dots,1)$; this is the line where all the inequalities $x_0 \geq x_j$,
$j = 1, \dots, n$ are equalities. Thus we can split $\calC_{n+1}$ as a
   sum of a line and lower dimensional cone
$$
\calC_{n+1} = \calL \oplus \calV_n,
$$
where $\calV_{n}=\calL^\perp\cap\calC_{n+1}$.
Notice that $V_j := e_0 - e_j$, $j = 1, \dots, n$ is orthogonal to
$(1,1,\dots,1)$, so $V_1, \dots, V_n$ provide a basis of
$\calL^\perp$. It is convenient to define
$$\calT_{n-1} := \calV_n \cap \Ss^{n-1},$$
where $\Ss^{n-1}$ is the unit sphere in $\calL^\perp = \Span \{ V_j\}$.
The domain $\calD_n$ is a
double cone over $\calT_{n-1}$. More precisely, let $N$ be one of the
intersection
points in $\Ss^n \cap \calL$ (there are two such points),
and let $(r,\theta)$ be polar coordinates in $\Ss^n$, centered at $N$. Then
$$
\calD_n = \{ (r,\theta) \mbox{ }|\mbox{ } \theta \in \calT_{n-1}, 0 \leq r
\leq \pi \},
$$
where we have identified $\Ss^{n-1}$ with the unit sphere in the tangent space
$T_N\Ss^{n}$.
In Figure \ref{fig-domains} we show $\calT_1$ and $\calD_2$.
\begin {figure}[h]
\begin {center}
\includegraphics[height=2in]{DandT1.eps}
\caption{This figure shows $\calT_1$ and $\calD_2$.}
\label {fig-domains}
\end {center}
\end {figure}

In later sections, we will use a generalization of this type of spherical
cone.
In general, let $\Omega$ be a domain in the equatorial $\Ss^{n-1}$
of $\Ss^n$, and
let $r_0 \in (0,\pi]$. Then we define the truncated cone
$$\calT\calC(\Omega,r_0) := \{(r,\theta) \mbox{ }|\mbox{ } \theta
\in \Omega, 0 \leq r \leq r_0 \}.$$
We abbreviate $\calT\calC(\Omega, \pi) = \calT\calC(\Omega)$. In this notation,
$\calD_n =
\calT\calC(\calT_{n-1}) = \calT\calC(\calT_{n-1},\pi)$.

The domain $\calT_{n-1}$ has symmetry. If we let
$$\calC_{n+1}^j = \{ (x_0, x_1,\dots x_n) \mbox{ }|\mbox{ } x_j \geq x_k,
j \neq k \},$$
then we see $\calC_{n+1} = \calC_{n+1}^0$ and $\calC_{n+1}^j$
are pairwise
congruent. Thus
$\calT_{n-1}$ is a face of the regular $(n+1)$-hedral tesselation of the
standard
$\Ss^{n-1}$ one obtains by connecting the vertices of a regular $(n+1)$-simplex
with great circle arcs. In particular, one can compute the diameter of
$\calT_{n-1}$
as the distance from a vertex to the center of the opposite face, which is
$$\delta(n-1) = \arccos (-\sqrt{\frac{n-1}{2n}}).$$
Moreover, the spherical angle of $\calT_{n-1}$ at a vertex is $\calT_{n-2}$.
So we
can construct a succession of comparison domains for $\calT_1,
\calT_2, \dots$ starting with $\calT_1$ and using
the coning construction described above. To this end, we let
$$\hat \calT_1 := [0,\frac{2\pi}{3}] = \calT_1, \qquad \hat \calT_n :=
\calT\calC(\hat \calT_{n-1}, \delta(n)).$$
By induction, $\calT_n \subset \hat \calT_n$, and so $\lambda_1(\calT_n) \geq
\lambda_1(\hat \calT_n)$.

\section {Separating variables} \label{separate-sec}

We  discuss two types of separation of variables in this section. The first
type is the separation of variables used in \cite{DB} to estimate expected exit
times of Brownian motion from Euclidean cones, and the second is the
separation
of variables one performs to estimate eigenvalues of a spherical domain with a
conical structure.

One needs the domains in question to have a certain amount of regularity
in order to separate variables.
In particular, these domains are piecewise $C^1$, satisfy an exterior
cone condition,
and are well approximated (in terms of Hausdorff distance) by finite
volume, $C^1$
domains. See the introduction to \cite{DB} for a precise statement. It is
straight-forward to verify that all of the domains we  consider
satisfy DeBlassie's
hypotheses.
We will refer
to these domains as {\em nice}.

\subsection{DeBlassie's separations of variables}

We first review DeBlassie's \cite{DB} argument. Consider the cone
$C\in\R^{n+1}$
over a domain $D \subset \Ss^n$:
$$
C = \{r\theta \mbox{ }|\mbox{ } r>0, \theta \in D \subset \Ss^{n} \}.
$$
Let $\tau_x$ be the time it takes for a Brownian particle to exit $C$, with
starting position $x$, and let $u(x,t) = \mathbb{P}(\tau_x >t)$ be the
probability that $\tau_x > t$. Then $u$ satisfies the heat equation
$$
\begin {array} {ll}
u_t = \frac{1}{2} \Delta u, & (x,t) \in C \times [0,\infty); \\
u(x,0) = 1, & x \in \bar C;\\
u(x,t) = 0, & (x,t) \in \del C \times (0,\infty).
\end {array}
$$
In polar coordinates $(r,\theta,t)$,  the PDE becomes
$$
2u_t = u_{rr} + \frac{n}{r} u_r + \frac{1}{r^2} \Delta_{\Ss^n} u.
$$
Moreover, the solution scales as $u(r,\theta,t) =
u(\beta r,\theta, \beta^2 t)$,
so we can separate variables and look for a solution of the form
$u= R(\xi) U(\theta)$, where $\xi = r^2/2t$. Then we get a
positive separation constant $\lambda_j(D)$, and
$$
\lambda_j (D) = \frac{\Delta_{\Ss^n} U}{U} = \frac{4\xi^2 \ddot R +
(4\xi^2 + 2n\xi)\dot R}{R}.
$$
Letting $R = \xi^a \rho(-\xi)$, this ODE becomes
$$
\xi \ddot \rho + (a + \frac{n+1}{2} - \xi) \dot \rho - a \rho = 0,
$$
which has solutions of the form $\rho(\xi) = f(a/2, a+(n+1)/2, \xi)$,
where $f$ is the confluent hypergeometric function
(see {\it e.g.}, \cite{lebedev}).
Using $f$ and restricting to $K\subset D$ compact and $T>0$, one obtains
(see \cite{DB}) an expansion for $u$ of the form
$$
u(x,t) = \sum_{j=1}^\infty B_j f \left ( a_j(n), 2a_j(n)+\frac{n+1}{2},
-\frac{|x|^2}{2} \right ) U_j \left ( \frac{x}{|x|} \right )
\left ( \frac{|x|^2}{2t} \right )^{-a_j(n)},
$$
where $U_j$ is the $j$th Dirichlet eigenfunction of $\Delta_{\Ss^n}$ on $D$ and
$$
2a_j(n) = \left [ \left (\frac{n-1}{2} \right ) ^2 + \lambda_j (D) \right
]^{1/2} - \frac{n-1}{2}.
$$
This sum converges uniformly in $K \times [T, \infty)$. 
The leading term in the expansion is $u(x,t) \sim B_1 U_1(x/|x|)
(|x|/(2t))^{-a_1(n)}$, which yields the conclusion of inequality
(\ref{eigen-est1}).

\subsection{Separating variables on the sphere}

Next we separate variables to relate the eigenvalues $\lambda_1(\calD_n)$ and
$\lambda_1(\calT_{n-1})$. First recall that we can write the
Laplacian for $\Ss^n$ in polar coordinates as
\begin {equation} \label {sep-var1}
\Delta_{\Ss^n} u = u_{rr} + (n-1) \cot r u_r + \csc^2 r \Delta_\theta u,
\end {equation}
where $\Delta_\theta$ is the Laplacian on the equatorial $\Ss^{n-1}$.
This lemma gives
us a dimension reduction. 

\begin {lemma}
Let $\Omega$ be a nice domain in an equatorial $\Ss^{n-1}$ with first
eigenvalue $\lambda = \lambda_1(\Omega)$, and let $\calD =
\calT\calC(\Omega)$ be the double cone over $\Omega$. Then
the first Dirichlet eigenvalues of $\calD$ and $\Omega$ are related by
\begin {equation} \label {eigen-rel1}
\lambda_1(\calD) = \lambda_1(\Omega) - \frac{n-2}{2} +
\sqrt{\frac{(n-2)^2}{4} + \lambda_1(\Omega)}.
\end {equation}
In particular, $\lambda_1(\calD) > 2n+2$ whenever $\lambda_1(\Omega) > 2n$.
\end {lemma}

{\bf Proof}: Set $u(r,\theta) = R(r) T(\theta)$, where $R(0) = 0 = R(\pi)$
and $T(\theta) = 0$ for $\theta
\in \del \Omega$. Then $u$ is an eigenfunction on $\calD$ with
eigenvalue $\mu$ precisely when
$$
T\ddot R + (n-1) \cot r\, T \dot R + \csc^2 r\, R \Delta_\theta T = - \mu TR.
$$
Separating variables with a positive separation constant $\lambda$ yields
$$
\frac{\sin^2 r \ddot R + (n-1) \sin r \,\cos r\, \dot R + \mu \sin^2 r R}{R}
= \lambda  = -\frac{\Delta_\theta T}{T}.
$$
Choosing $T$ to be the first eigenfunction of $\Omega$, we obtain the ODE
\begin {equation} \label {sep-var2}
\sin^2 r \ddot R + (n-1)\sin r \cos r \dot R + (\mu \sin^2 r - \lambda)R = 0,
\end {equation}
which has regular singular points at $r = 0,\pi$. If we try a
solution of the form $R = \sin^m r$, for
some power $m$, we find
\begin {eqnarray*}
0 & = & m(m-1)\sin^m r \cos^2 r- m \sin^{m+2} r + m(n-1)\sin^m r\cos^2 r +
(\mu \sin^2 r - \lambda) \sin^m r \\
& = &\sin^m r [(m^2 + m(n-2) - \lambda) \cos^2 r + (\mu - m -
\lambda)\sin^2 r ].
\end {eqnarray*}
Both coefficients must vanish, so we have $\mu = m+\lambda$ and
$\lambda = m^2 + m(n-2)$.
Solving for $m$, we find
\begin {equation} \label {rel-eigenval}
m = \frac{2-n}{2} + \sqrt{\frac{(2-n)^2}{4} + \lambda }.
\end {equation}
Next, observe that if $\lambda_1(\Omega) = 2n$ then $\lambda_1(\calD) =
2n+2$. Finally, the formula for
$\lambda_1(\calD)$ is monotone increasing in $\lambda_1(\Omega)$,
and so $\lambda_1(\calD) > 2n+2$
whenever $\lambda_1(\Omega) > 2n$.
\hfill $\square$

\begin {rmk}
A second solution to equation (\ref{sep-var2}) has the form $\sin^m r
\cos r$, where $m$ is again given
by equation (\ref{rel-eigenval}) but
$$
\mu' = \lambda + 3m + n.
$$
This eigenfunction vanishes on $\{ \pi/2 \} \times \Omega$, so it
corresponds to a higher eigenvalue.
\end {rmk}

At this point, we can prove that the expected capture time for $n=1,2,3$
predators
is infinite. To prove the expected capture time is infinite, by
inequality (\ref{eigen-est1}) we need to show $\lambda_1(\calD_n)
\leq 2n+2$, or,
equivalently, that $\lambda_1(\calT_{n-1}) \leq 2n$. In the case of
$n=1$, we have
$\calD_1 = [-3\pi/4, \pi/4]$, and so $\lambda_1 (\calD_1) = 1 < 4$.
In the case
$n=2$, we have $\calT_1 = [0, 2\pi/3]$, and so $\lambda_1(\calT_1) =
9/4 < 4$. We
cannot compute $\lambda_1(\calT_2)$ so easily, but we can find a test
function to
show that $\lambda_1(\calT_2) < 6$. We shall find
$\lambda_1(\calT_2)$ numerically in
Section~\ref{numer-sec}. Recall the Rayleigh characterization  of the
first eigenvalue of a domain $\Omega$:
$$
\lambda_1(\Omega) = \inf_{f \in H^1_0(\Omega), f \not \equiv 0}
\left ( \frac{\int_\Omega |df|^2} {\int_\Omega f^2} \right ).
$$
To show that $\lambda_1(\calT_2) < 6$, it suffices to find $f_0
\in H^1_0(\calT_2)$
so that $\int_{\calT_2} |df_0|^2 /\int_{\calT_2} f_0^2 < 6$. Let
$$
f_0(x) = \sin (\dist(x,\del \calT_2)),
$$
and observe that $|df_0|^2 = 1 - f_0^2$ off the set of focal points of
$\del \calT_2$, which is a set of measure zero. A computation shows
$$
\lambda_1(\calT_2) \leq
\frac{\int_{\calT_2} |df_0|^2}{\int_{\calT_2} f_0^2} = \frac{2\pi + \sqrt{3}}
{\pi - \sqrt{3}} < 6.
$$

One can generalize the eigenvalue relationship (\ref{eigen-rel1}) to
spherical cones of the form $\calT\calC(\Omega,r_0)$, for $0 < r_0 < \pi$,
using the hypergeometric function
\begin {equation} \label{confluent}
f(\alpha, \beta, \gamma, z) = 1 + \frac{\alpha}{\beta} \frac{z}{\gamma} +
\frac{\alpha(\alpha + 1)}{\beta(\beta+1)} \frac{z^2}{\gamma(\gamma+1)} +
\frac{\alpha(\alpha+1)(\alpha+2)}{\beta(\beta+1)(\beta+2)}
\frac{z^3}{\gamma(\gamma+1)(\gamma+2)} + \dots
\end {equation}
Later we will use the following lemma to
relate the eigenvalues of $\hat \calT_{n-1}$ and $\hat \calT_n$, which
we will then estimate to bound the asymptotic decay rates of probabalistic
exit times.
\begin {lemma}
Let $\Omega$ be a nice domain in a equatorial $\Ss^{n-1} \subset \Ss^n$ with
eigenvalue $\lambda = \lambda_1(\Omega)$. Then the
first Dirichlet eigenvalues of $\Omega$ and $\calT\calC(\Omega,r_0)$ are
related by
\begin{equation} \label{eigen-rel2}
\lambda_1(\calT\calC(\Omega,r_0)) = \mu = \mu(n, \lambda, r_0),
\end {equation}
where $\mu$ is the first Dirichlet eigenvalue of the ODE (\ref{sep-var2}) on
$[0,r_0]$.  If $r_0 \geq \pi/2$ then $\mu$ is the unique zero of
$f(\alpha_1, \beta_1,
\gamma_1, (1/2)(1- \cos r_0))$ in $(m+ \lambda, 3m + \lambda + n)$, where
$m$ is defined in equation (\ref{rel-eigenval}) and
\begin {eqnarray*}
\alpha_1, \beta_1 & = & \frac{ 1 + \sqrt{(n-2)^2 + 4\lambda}
\pm \sqrt{(n-1)^2 + 4\mu}}2 \\
\gamma_1 & = & \frac{2 + \sqrt{(n-2)^2 + 4\lambda}}{2}.
\end {eqnarray*}
If $r_0 \leq \pi/2$ then $\mu$ is the unique zero of $f(\alpha_2, \beta_2,
\gamma_2, (1/2)(1-\cos r_0))$ in $(0,n)$, where
$$
\alpha_2,\beta_2 = \frac{n-1 \pm \sqrt{(n-1)^2 + 4\mu}}{2}, \qquad
\gamma_2 = \frac{n}{2},
$$
with $f$ defined by equation (\ref{confluent}).
\end {lemma}

{\bf Proof}: We separate variables and look for a solution of the form
$R(r) = \sin^m r\, u(r)$, with $u(r) \neq 0$ on $[0,r_0)$, but $u(r_0) = 0$.
Then equation (\ref{sep-var2}) becomes
$$
0 = \sin^{m+2} r\, \ddot u + (2m+n-1) \sin^{m+1} r \,\cos r\, \dot u +
[m(m+n-2) \cos^2 r + (\mu - m)\sin^2 r - \lambda] \sin^m r\, u.
$$
Now let $u(r) = y(x)$, where $x = (1/2)(1 - \cos r)$, which transforms
the ODE above into
$$
x(1-x)y'' + (m+ \frac{1}{2} n - (2m+n)x) y' - (\lambda + m - \mu)y = 0.
$$
The solution to this ODE is the hypergeometric function $y(x) =
f(\alpha, \beta, \gamma, x)$, with
\begin {eqnarray*}
\alpha, \beta & = & \frac{2m + n - 1 \pm \sqrt{(2m+n-1)^2 - 4\lambda
-4m+4\mu}}{2} \\
\gamma & = & \frac{2m+n}{2}.
\end {eqnarray*}
The lemma follows from taking $R(r;n,\lambda, r_0) =
\sin^m r\, f(\alpha,\beta,\gamma,(1/2)(1-\cos r))$, where we choose
$\mu$ so that $R(r_0; n, \lambda, r_0) = 0$.
\hfill $\square$

One can use this lemma to compute $\lambda_1(\hat \calT_{n-1})$
iteratively. In this case we start with $\lambda_1(\hat \calT_1) = 9/4$
and apply equation (\ref{eigen-rel2}).
For the reader's reference, we include a table of $\lambda_1(\hat \calT_{n-1})$
and the corresponding lower bound for $a(n)$ (see equation (\ref{decay-exp})).
In this
table we use equation (\ref{eigen-rel1}) and the fact that $\calD_n \subset
\calT\calC(\hat \calT_{n-1})$ so $\lambda_1(\calD_n) \ge
\lambda_1(\calT\calC(\hat \calT_{n-1}))$.
$$
\begin {array}{cll}
\underline n & \underline {\lambda_1(\hat \calT_{n-1})} &
\underline {\mbox{lower bound for }a(n)} \\
2 & 2.25 &\ \enspace .75000000\\
3 & 5.00463581 &\ \enspace .89614957\\
4 & 7.884040724 &\ \enspace .99030540\\
5 & 10.77018488 &\ 1.05417466\\
6 & 13.6203196 &\ 1.09882819
\end {array}
$$
  From this table, one can see that $\lambda_1(\hat \calT_3)$ is
slightly less than $8$, so $\hat \calT_3$ comes close to, but is
not quite, an effective comparison domain for proving
Theorem \ref{main-thm}.

\section {The eigenvalue estimate} \label{eigen-est-sec}

In this section, we complete the proof of Theorem \ref{main-thm}. First
observe
that it suffices to show $\lambda_1(\calT_3) > 8$. In fact, it
suffices to find an
effective lower bound for $\lambda_1(\calT_2)$ using the following scheme.
We define
$\lambda_{cr}$ by
$$
8 = \mu (3,\lambda_{cr}, \delta(3)), \qquad
\lambda_{cr} \simeq 5.101267527.
$$
   If we construct a domain
$\calG_2 \subset \Ss^2$  such that
$\calT_2 \subset \calG_2$ and $\lambda_1(\calG_2)>\lambda_{cr}$, then by domain
monotonicity and equation (\ref{eigen-rel2}),
$$
\lambda_1(\calT_3)\ge
\lambda_1(\calT\calC(\calT_2,\delta(3))))\ge\lambda_1(\calT\calC(\calG_2,\delta(3)))>8,
$$
completing the proof of Theorem \ref{main-thm}.
As the last step in our proof, we construct  a domain $\calG_2 \subset
\Ss^2$ as a
perturbation of $\calT_2$, such that $\calT_2\subset\calG_2$  and
$\lambda_1(\calG_2)=5.102 >
\lambda_{cr}$.

\begin {figure}[h]
\begin {center}
\includegraphics[width=4.1in]{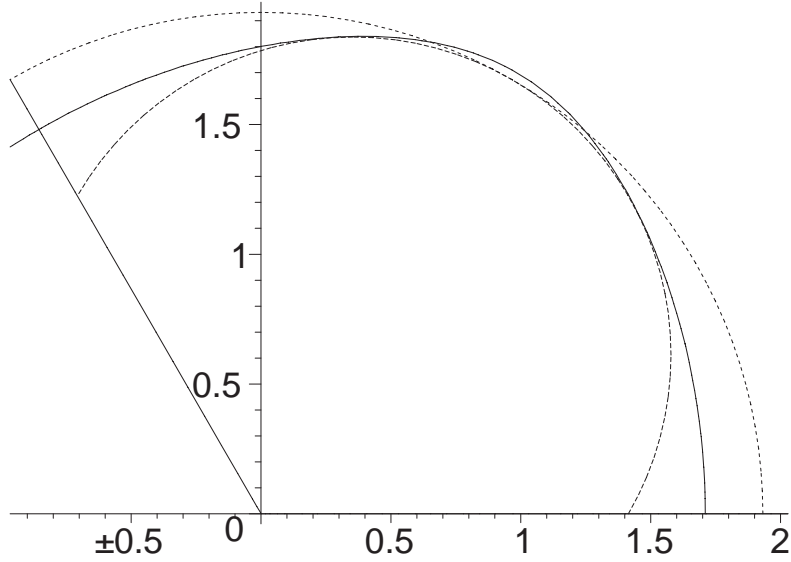}
\caption{The inner dashed curve is $\calT_2$, solid curve is $\calG_2$,
and the outer dashed curve is $\hat \calT_2$. We generated this
figure using the
computer program MAPLE.}
\label {fig-arcs}
\end {center}
\end {figure}

Figure \ref{fig-arcs} shows sketches of  the domains $\calT_2$, $\calG_2$ and
$\hat \calT_2$, after stereographic projection to the plane.

Previous work of Rayleigh
\cite {R} and P\'olya--Szeg\"o \cite {PS} motivates us to consider this
type of
domain perturbation. They studied the eigenvalue of a planar domain
which has the
form $\{(r,\theta) \mbox{ }|\mbox{ } 0 \leq r \leq c + \epsilon f(\theta) \}$
in polar coordinates, for some small $\epsilon > 0$, and gave an expression for
the eigenvalue $\lambda_1$ in terms of $\epsilon f$. In our case,
we fix $\lambda = \lambda_1$ and find a domain $\calG_2$ with
$\lambda$ as its first eigenvalue.

Suppose the functions $R(r)$ and $\Theta(\theta)$ satisfy
$$
\begin {array}{rl}
\Theta'' + \lambda \Theta = 0, & 0 \leq \theta \leq \frac{2\pi}{3}; \\
\Theta(0) = 0,\\  \Theta(2\pi/3) = 0, \\
\sin^2 r\, \ddot R + \sin r \,\cos r \,\dot R + (\mu \sin^2 r -
\lambda) R = 0, &
0 \leq r < \pi; \\
R(0) = 0. & \end {array} $$
Then $u(r,\theta) = R(r) \Theta(\theta)$ is the first eigenfunction of
$\Delta_{\Ss^2}$ on its nodal domain $\calG_2$. By construction,
$\Delta_{\Ss^2} u
+ \mu u = 0$, so $u$ is an eigenfunction. Also, $u$ does not change
sign on its nodal domain, so it must be the first eigenfunction. In polar
coordinates, $u(r,0) = 0 = u(r, 2\pi/3)$. Let $m =
\sqrt{\lambda}$ and
set $R(r) = \sin^m r u(r)$, so that equation (\ref{sep-var2}) becomes
$$
\sin^2(r) \,\ddot u + (1 + 2m) \sin (r)\, \cos(r) \,\dot u +
(\mu - m - \lambda) u = 0.
$$
Next we take $\lambda = 9l^2/4$, corresponding to the $l^{th}$ mode
of the interval $[0, 2\pi/3]$, and write the solution in terms of the
hypergeometric function:
$$u_l (r) = f(3l/2 + .5 \pm \sqrt{1/4 + \mu}, 1 + 3l/2, (1-\cos r)/2).$$
Finally, we take the $\mu = 5.102$ superposition of the $l=1,3$ modes
to define
$$\sin^{3/2} (r) \sin(3\theta/2) H(r, \theta)  := (\sin r)^{3/2} u_1(r)
\sin (3\theta/2) - .0003 (\sin r)^{9/2} u_3(r) \sin (9\theta/4)$$
and let $\calG_2$ be the nodal domain of $\sin^{3/2}(r) \sin(3\theta/2) H$.
By construction, $\lambda_1(\calG_2) = 5.102 > \lambda_{cr}$.

The domains $\calT_2$, $\calG_2$ and $\hat \calT_2$ all have three boundary
curves, and two of these boundary curves for each domain lie along the
great circle arcs $\theta= 0, 2\pi/3$. For convenience, we convert to a
planar domain using stereographic projection, with the south pole
corresponding to $r=0$ in polar coordinates on $\Ss^2$. The relationship
between polar coordinates $(r,\theta)$ on $\Ss^2$ and polar coordinates
$(\rho, \theta)$ on the plane is $\rho = \tan(r/2)$. Under stereographic
projection, the great circle arcs $\theta = 0, 2\pi/3$ correspond to
rays emanating from the origin at angles $0,2\pi/3$. The third boundary
curve of $\calT_2$ is given by the circular arc
$$
(\beta(\theta) \cos\theta - \frac 1{\sqrt{8}})^2 + (\beta(\theta) - \sqrt{\frac
32})^2  = \frac{3}{2},
$$
which we can rewrite as
$$
\rho = \beta(\theta) = \frac{\sqrt{2}\cos(\theta - \pi/3) +
\sqrt{2 \cos^2 (\theta - \pi/3) + 4}}{2}.
$$
Thus showing $\calT_2 \subset \calG_2$ is equivalent to showing
$H(r,\theta)>0$ along the arc given by $(\rho = \beta(\theta), \theta)$ for
$0 \leq \theta \leq 2\pi/3$.

One can see from Figure \ref{fig-arcs}
that $H>0$ along this third boundary
component of $\calT_2$. One strategy for a rigorous proof that $H>0$ is
the following. We first evaluate $H$ at a point $\theta_0$
on the curve, checking $H>0$ at $(\rho = \beta(\theta_0), \theta_0)$,
and bound the derivative $H_\theta$ on an interval containing
$\theta_0$. Our bound $H_\theta \geq M$ gives us a lower bound
$H > H(\theta_0) - M (\theta - \theta_0)$. Thus $H> 0$ on a possibly smaller
neighborhood of $\theta_0$. We then repeat this procedure with each endpoint
of this (smaller) interval. One can simplify the computations by observing
$\calT_2$ and $\calG_2$ are symmetric about the ray $\theta = \pi/3$. Then it
suffices to evaluate $H$ and its derivative at $\theta = 0, 1/2, 2/3, 2\pi/9, 
\pi/3$.
\hfill $\square$.

Given our lower bound $\lambda_1(\calT_2) \geq \lambda_1(\calG_2) = 5.102$,
one can use equation (\ref{eigen-rel2}) to show $\lambda_1(\calT_3) \geq
8.00087815$. This in turn gives $\lambda_1(\calD_4) \geq 10.001024501$ and
$a(4) \geq 1.00007318$.
In the same way, the lower bound $\lambda_1(\calT_2) \geq 5.102$ gives
$a(3) \geq .90671950$.


\section {A numerical computation of the first eigenvalue} \label{numer-sec}

In this section we describe a numerical computation approximating
the eigenvalue $\lambda_1(\calT_2)$, which relies on Stenger's
sinc-Galerkin scheme \cite{St}. We will show $\lambda_1(\calT_2) \approx
5.159\dots$, and so $a(3) \approx .9128\dots$. Also, using equation
(\ref{eigen-rel2}), we have $\lambda_1(\calT_3) \stackrel
{\textstyle >}{\sim} 8.0691$ and so
$a(4) \stackrel{\textstyle >}{\sim} 1.0057$.

Given $h>0$ and a positive integer $k$, one defines the $k$th cardinal sinc
function of stepsize $h$ as $S(h,k)(hk)=1$ and if $z\ne hk$,
\begin {equation} \label{sinc-eq}
S(h,k)(z) := \frac{\sin(\pi(z-hk)/h)}{\pi(z-hk)/h}.\end {equation}
Following Stenger\cite{St}, for $z=x+iy$ we shall let the basis functions be
$\phi_{jk}(z)=\alpha_j(x)\times \beta_k(y)$, where
\begin {equation} \label{sinc-basis1}
\alpha_j(x)=S(j,h)\circ \ln (\frac{x}{\frac\pi 2-x}),\quad
\alpha_{n+1}(x)=\sin^2(x)-\sum_{l=-n}^n \sin^2(x_l) \alpha_l(x)
\end {equation}
for $j=-n\ldots,n$ with sinc points $x_l=\dfrac{\pi e^{hl}}{2(1+e^{hl})}$ and
\begin {equation} \label{sinc-basis2}
\beta_k(x)=
S(k,h)\circ \log\sinh  y,\quad
\beta_{n+1}(y)=\sech(y)-\sum_{l=-n}^n \sech(y_l) \beta_l(y),
\end{equation}
for
$k=-n,\ldots,n$ with sinc points $y_l=\sinh^{-1}(e^{hl})$.
We let $x_{n+1}=\pi/2$  and
$y_{n+1}=0$. The dimension  of the space of functions $\phi_{jk}$
is $m=(2n+2)^2$. Also $\alpha_{n+1}(\pi/2) =1$ and
$\beta_{n+1}(0) = 1$, but both are zero at the other sinc points because
their interior sinc
approximations are subtracted. The functions $\phi_{jk} = \alpha_j \times
\beta_k$ form a basis of sinc functions when we are working on a
rectangular strip.

The most important property of sinc functions is that they are
well-suited to
approximating integrals in a strip. In particular,
one can increase accuracy in a numerical computation by increasing the number
of evaluation points, without  recomputing grids (as in finite 
elements). In order to take
advantage of this simplicity,
we will
conformally map a subdomain of $\calT_2$ to a half-infinite strip.

\begin {figure}[h]
\begin {center}
\includegraphics [height = 1.5in]{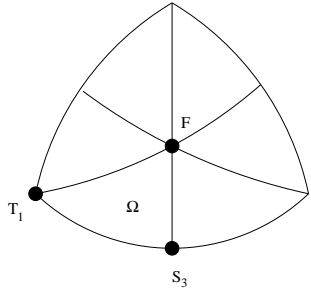}
\caption {The three marked points are $F$, $T_1$ and $S_3$.}
\label{fig-triangles}
\end {center}
\end {figure}

Recall that $\calT_2$ is an equilateral triangle on $\Ss^2$ with all its
interior angles equal to $2\pi/3$. Let $T_1,T_2,T_3$ be the vertices of
$\calT_2$, let $S_j$ be the midpoint of the side opposite $T_j$, and let
$F$ be the center of mass of $\calT_2$. Observe that $\calT_2$ is invariant
under reflection through the lines $FT_j$, which divide $\calT_2$ into
six congruent subtriangles (see Figure \ref{fig-triangles}).
Thus the first eigenfunction is invariant
under these reflections, and we can recover it by restricting to the
smaller triangle $\Omega$, which has vertices $F,T_1,S_3$. The triangle
$\Omega$ has a right angle at the vertex $S_3$, and an angle of $\pi/3$ at
the vertices $F$ and $T_1$. The first Dirichlet eigenfunction of $\calT_2$,
restricted to $\Omega$,  will have Dirichlet data on the edge $T_1 S_3$ and
Neumann data on the edges $F T_1$ and $F S_3$.

We first transform $\Omega$ to
a planar domain (which we again denote as $\Omega$) using stereographic
project, sending $S_3$ to $0$. This
transformation sends $F$ to $\frac{1}{2}(\sqrt{6} - \sqrt{2})$ and $T_1$
to $\frac{i}{2}(\sqrt{6} - \sqrt{2})$. Now
$\Omega$ is bounded by the two straight line
segments joining $0$ to $\frac{1}{2}(\sqrt{6} - \sqrt{2})$ and $0$ to
$\frac{i}{2}(\sqrt{6} - \sqrt{2})$ and the circular arc joining $\frac{1}{2}
(\sqrt{6}- \sqrt{2})$ to $\frac{i}{2}(\sqrt{6} - \sqrt{2})$ which makes an
angle of $\pi/3$ with the axes. Next we find a conformal transformation $f$
which takes $\Omega$ to the half-infinite strip
$$D = \{z \in \C \mbox{ }|\mbox{ } 0 < \Re z < \pi/2, 0 < \Im z\}$$
such that $f(0) = 0$, $f(\frac{i}{2}(\sqrt{6} - \sqrt{2})) = \pi/2$, and
$f(\frac{1}{2}(\sqrt{6} - \sqrt{2})) = \infty$. Under this transformation,
the first eigenfunction $u$ satisfies
\begin{eqnarray*}
\Lap^* u +\lambda u&=&0, \qquad \text{if $z\in D$,}\\
u&=&0,\qquad\text{if $\real z=0$ and $\imaginary z>0$,}\\
\frac {\partial u}{\partial n} &=&0,\qquad \text{if $0<\real  z<\frac 
\pi 2$ and
$\imaginary z =0$, or if $\real z=\frac \pi 2$ and
$\imaginary z >0$;}
\end{eqnarray*}
where $\Delta^*=f^*\Delta=\Delta_z$ is the pulled back Laplacian.

There are formul\ae\  for conformally mapping  domains bounded by finitely many
circular arcs \cite{N} generalizing the Schwarz-Christoffel formula. Following
Forsythe \cite{F}, the Schwarz triangle mapping $z\in D$ or from
$\sin^2 z\in\{\Im z > 0\}$ of the upper halfplane to
$w\in\Omega$ is given by
$$
\cos^2 z=\frac{(w^4+2\sqrt{3}\, w^2-1)^3}{(w^4-2\sqrt{3}\, w^2-1)^3}.
$$
   The group generated by   reflections of
$\Omega$ along its edges in fact tiles the sphere. The inverse function
$w\mapsto z^2$ extends to  a single sheeted cover and is invariant under
the symmetry group,
thus is a rational function of the plane.
Thus we may compute $f$. Writing $g(z)=\cos^{2/3}z$,
$$
f(z)=
\sqrt{
\frac{1-g}{\sqrt 3(1+g) +2\sqrt{1+g+g^2}}
}
$$
Pulling back under $w=f(z)$, the conformal weight takes the form
$$
\frac{4\left| \dfrac{df}{dz}\right|^2}{( 1 + |f|^2)^2}= \frac{\dfrac 43
\left| \sqrt 3
(1+g)+2\sqrt{1 + g + g^2}\right| }{|g|(\left| \sqrt 3 (1+g)+2\sqrt{1 + g +
g^2}\right|+|1-g| )^2}.
$$
The branch cuts for the square and cube roots may be taken above the
negative real axis. Thus $g(D)$ lies in the fourth
quadrant so the denominator in $f$ is nonvanishing.

Let $G(z,z')$ denote the Green's function for the problem on $D$
\begin{eqnarray*}
\Lap^* u &=&\, g, \qquad \text{if $z\in D$,}\\
u&=&0,\qquad  \text{if $\real z=0$ and $\imaginary z>0$,}\\
\frac {\partial u}{\partial n} &=&0,\qquad\text{if $0<\real z < \pi/2$
and $\imaginary z=0$ or if
$\real z = \pi/2$ and $\imaginary z>0$.}
\end{eqnarray*}

The Green's function may be found by the method of images. Thus, denoting
$w=\sin z=x+iy$,
$w^*=-x+yi$ and $\omega=\sin\zeta$, we get $\overline{w^*}=(\bar w)^*$.
The Green's function is
$$
G(z;\zeta)=\frac 1{2\pi}\Bigl({\ln|w-\omega|-\ln |w^\ast-\omega| +
\ln|\bar w-\omega|
-\ln|\bar w^\ast - \omega |}\Bigr).
$$
  Pulling back by $f$, the eigenvalue problem for the triangle may be
restated as finding an eigenvalue for the integral operator
\begin{equation}\label{tag8}
\frac 1\lambda u(z) = -4\int_D  \frac{G(z,z')\, |df(z')|^2\, u(z')\,
dz'}{(1+|f(z')|^2)^2}=:\mathcal A u(z)
\end{equation}

  The key point is to notice that the operator has logarithmic and algebraic
singularities at the points $z'=0,\pi/2$ and $z= z'$. Thus we need to 
handle these
singularities. The solution is zero along the imaginary axis, but free
along the other
two sides, which may be extended to functions to the plane which have odd
reflection
symmetry along the imaginary axis and even reflection symmetry along the
other sides.
We shall approximate
$u(z)$ in an
$m$-dimensional space
  $X_m$ with the same symmetries. Also, noticing that the
eigenfunction on $\mathcal
T_2$ behaves like $1-\dist^2(z,T_2)$ at the vertex $F$, we actually have $u\in
\Lip(\bar D)$ and $u$ decays algebraically at $\infty$. Choosing a
  basis
$\{
\phi_1,\ldots,\phi_m\}$ of
$X_m$, we shall compute the matrix of the transformation
  $A_{\ell k}=\mathcal P_\ell\mathcal A \phi_k$, whose largest eigenvalue
$\mu_m\to 1/\lambda$ as
$m\to\infty$ and which is an upper bound $\lambda \le 1/\mu_m$ \cite{St}.
The integral
operator shall be computed numerically via sinc quadrature, which can
handle such mild
singularities.

We make the approximation $u(z)\approx \mathcal P u(z)= b^{jk}\phi_{jk}(z)$
where the sums are over $j=-n\ldots,n+1$ and
$k=0,\ldots,n+1$. Here $\phi_{jk} = \alpha_j \times \beta_k$, where
$\alpha_j$ and $\beta_k$ are defined in equations (\ref{sinc-basis1}) and
(\ref{sinc-basis2}), respectively.
We define $b^{jk}=\mathcal P_{jk} u = u(x_j+iy_k)$.
Thus the approximation $\mathcal P  u$ is a collocation, as it  equals
$u$ at the sinc points.
Thus, the matrix is approximated by
\begin{eqnarray*}
A_{jk,pq}&=&\int_D G(x_j,y_k,\xi,\eta)\phi_{pq}(\xi,\eta) \Psi(\xi,\eta)\,
d\xi\,d\eta\\
&\approx &  \sum _{\iota,\kappa} v_\iota w_\kappa
G(x_j,y_k,x_\iota,y_\kappa)\phi_{pq}(x_\iota,y_\kappa)\,
\Psi(x_\iota,y_\kappa)
\end{eqnarray*}
where
$
\Psi(\xi,\eta)={4 {(1+|f(\xi+i\eta)|^2)^{-2}} |df(\xi+i\eta)|^2}
$.
The approximating sum is carried over $4m$ terms, corresponding to sinc
quadratures in the four regions bounded by
singularities ({\it e.g.\,} in case $0<x_j<\pi/2$ and $0<y_k$):
$
D_I=\{ \xi+i\eta: 0<\xi<x_j,\ 0<\eta<y_k\}$, $
D_{II}=\{ \xi+i\eta: x_j<\xi<\pi/2,\ 0<\eta<y_k\}$, $
D_{III}=\{ \xi+i\eta: 0<\xi<x_j,\ y_k<\eta \}$ and $
D_{IV}=\{ \xi+i\eta:  x_j<\xi,\ y_k<\eta\}$ and where $v_\iota$ and
$w_\kappa$ are
corresponding weights and $x_\iota$ and
$y_\kappa$ are
corresponding sinc points.

The computation shows that $\lambda_1(\mathcal T_2)\approx 5.159\ldots $
which is above
the critical value. This is  computational evidence that the exit 
time has finite
expectation. We present a table of a few the computed  eigenvalues for
approximations
in spaces of given dimension, coming from the Sinc-Galerkin collocation scheme
described.
$$
\begin{array}{cc}
\underline{\text{ Dimension}\phantom{g}} &
\underline{\text{ Eigenvalue Estimate }}\\
      16          &           5.948293885960918 \\
      100         &           5.262319373675790 \\
     1024         &            5.153693139833067 \\
     2116         &           5.158585939808193 \\
     2304         &            5.158832705984016 \\
    2500         &            5.158849530710926 \\
      2704         &           5.158968860560663 \\
\end{array}
$$


\begin {thebibliography}{99}

\bibitem [BG]{BG} M. Bramson and D. Griffeath. {\em Capture problems
for coupled random walks}.
in Random Walks, Brownian Motion and Interacting Particle Systems,
ed. R. Durrett and H. Kesten,
  Birkh\"auser , 1991.

\bibitem [DB]{DB} R. D. DeBlassie. {\em Exit times from cones in $\R^n$ of
Brownian motion}. Prob.
Theory and Rel. Fields. 74:1--29, 1987.

\bibitem [F]{F} A. R. Forsythe. {\em Theory of Functions of a Complex
Variable, 3rd ed.}
Dover, 1965, 698--700 (orignially published by Cambridge University Press,
Cambridge, 1918).

\bibitem[L]{lebedev} N. Lebedev, {\em Special Functions \& Their 
Applications}. Dover, 1972,
260--262.  (Originally published by Prentice-Hall, Inc., 1965.)


\bibitem [LS]{LS} W. Li and Q.-M. Shao. {\em Capture time of Brownian
pursuits}. Prob. Theory and Rel. Fields. 121:30--48, 2001.

\bibitem [N]{N} Z. Nehari. {\em Conformal Mapping}. Dover, 1975
(originally published by Mc Graw-Hill Book
Co., Inc., New York, 1952).


\bibitem [PS] {PS} G. P\'olya \&\  G. Szeg\"o. {\em Isoperimetric
Inequalities in Mathematical physics}.
Princeton University Press, 1951, 29--30.

\bibitem [R]{R} Lord Rayleigh (J. W. Strutt) {\em The Theory of Sound
(2nd. ed.)}
Dover, 1945, 336--342 (Originally published by Mcmillan Co., 1894--1896).

\bibitem [St]{St}  F. Stenger {\em Numerical Methods Based oN Sinc
and Analytic Functions}
Springer-Verlag, 1993.


\end {thebibliography}
\end{document}